\documentclass[a4paper, 10pt, conference]{ieeeconf}      

\IEEEoverridecommandlockouts                              
\usepackage{graphicx} 
\usepackage{amsmath} 
\usepackage{amssymb}  
\usepackage{color,soul}
\usepackage[dvipsnames]{xcolor}
\usepackage{algorithm}
\usepackage{algpseudocode}
\algtext*{EndFor}
\algtext*{EndIf}

\usepackage{subcaption} 
\usepackage{hyperref}
\usepackage{enumerate}

\newtheorem{defn}{\bf Definition}
\newtheorem{prop}{\bf Proposition}[section]
\newtheorem{lemma}{\bf Lemma}[section]

\newtheorem{thm}{\bf Theorem}[section]
\newtheorem{corol}{\bf Corollary}[section]
\newtheorem{lemmaappendix}{\bf Lemma}[subsection]

\renewenvironment{proof}{{\bfseries Proof:}}{\cqfd}

\providecommand{\com[2]}{\begin{tt}[#1: #2]\end{tt}}

\providecommand{\abs[1]}{\left|#1\right|}

\definecolor{mygreen}{RGB}{204,255,153}
\definecolor{myblue}{RGB}{153,204,255}
\definecolor{myorange}{RGB}{255,204,153}
\newcommand{\new}[1]{#1} 
\newcommand{\neww}[1]{#1} 
\newcommand{\newg}[1]{#1} 
\newcommand{\newb}[1]{#1} 
\newcommand{\newo}[1]{#1} 

\DeclareMathOperator{\newdiff}{d} 
\newcommand{\diff}{\newdiff\!}
\newcommand{\fpart}[2]{\frac{\partial #1}{\partial #2}}
\providecommand{\prt}[1]{\ensuremath{\left( #1 \right)}}
\newcommand{\cqfd}{\hfill \rule{2mm}{2mm}\medbreak\indent}
\newcommand{\evec}{\mathbf{e}}

\DeclareMathOperator{\newrank}{rank} 
\newcommand{\rank}{\newrank\,}
\DeclareMathOperator{\newvec}{vec} 
\newcommand{\myvec}{\newvec}
\DeclareMathOperator{\newmat}{mat} 
\newcommand{\mymat}{\newmat}

\newcommand{\mytexttt}[1]{{\fontfamily{lmtt}\selectfont#1}}
\newcommand{\myfct}{h}

\title{\LARGE \bf
Local Network Identifiability with Partial Excitation and Measurement
}

\author{Antoine Legat and Julien M. Hendrickx
\thanks{\newg{A. Legat and J. M. Hendrickx are with ICTEAM Institute, UCLouvain, Belgium.}
\newb{Work supported
by the ``RevealFlight'' ARC at UCLouvain, and by the Incentive Grant
for Scientific Research (MIS) ``Learning from Pairwise Data'' of the
F.R.S.-FNRS.} {\tt\small antoine.legat@uclouvain.be, julien.hendrickx@uclouvain.be.}}%
}

\begin{document}

	\maketitle
	\thispagestyle{empty}
	\pagestyle{empty}

	\begin{abstract}

This work focuses on the identifiability of dynamical networks with partial excitation and measurement: a set of nodes are interconnected by unknown transfer functions according to a known topology, some nodes are subject to external excitation, and some nodes are measured. The goal is to determine which transfer functions in the network can be recovered based on the input-output data collected from the excited and measured nodes. 

We propose a local version of network identifiability, representing the ability to recover transfer functions which are approximately known, or to recover them up to a discrete ambiguity. We show that local identifiability is a generic property, establish a necessary and sufficient condition in terms of matrix generic ranks, and exploit this condition to develop an algorithm determining, with probability 1, which transfer functions are locally identifiable. Our implementation presents the results graphically, and is publicly available.
	\end{abstract}

	\medskip
	
	\section{INTRODUCTION}
	
	\medskip

This paper addresses the identifiability of dynamical networks in which the node signals are connected by causal linear time-invariant transfer functions, and can be excited and/or measured. Such network can be modeled as a directed graph where each edge carries a transfer function, and known excitations and measurements are applied at \newg{certain} nodes.

\medskip

		We consider the \newb{identifiability} of a network matrix $G(z)$, where the network is made up of $n$ nodes, \newg{with node signals $w(t)$, external excitation signals $r(t)$, measured nodes $y(t)$ and noise $v_1(t), v_2(t)$ related to each other by:}
		\begin{align} \label{eq:networkModel} \begin{split}
			w(t) &= G(z)w(t) + Br(t) + v_1(t)\\
			y(t) &= Cw(t) + v_2(t),
		\end{split} \end{align}
		where matrices $B$ and $C$ are binary selections defining respectively which nodes are excited and measured, forming sets $\mathcal{B}$ and $\mathcal{C}$ respectively. Matrix $B$ is full column rank and each column contains one 1 and $n-1$ zeros. Matrix $C$ is full row rank and each row contains one 1 and $n-1$ zeros.
		The nonzero entries of the network matrix $G(z)$ define the network topology, and are the transfer functions to identify, forming the set of edges $E$.
		
		\medskip
		
		We assume that \emph{the input-output relations between the excitations $r$ and measures $y$ have been identified}, and that the network topology is known. From this knowledge, we aim at recovering the nonzero entries of $G(z)$.
		
		\medskip
	
	The model \eqref{eq:networkModel} has recently been the object of a significant research effort.
	If the whole network is to be recovered, the notion of network identifiability is used, as developed in \cite{weerts2015identifiability}. If one is interested in identifying a single module, topological conditions are derived in \cite{weerts2018single}, \cite{gevers2018practical}. In this paper, we do not consider the impact of noise signals $v_1,v_2$, but studying the influence of rank-reduced or correlated noise yields interesting results on identifiability \cite{weerts2018prediction}, \cite{van2019local}.
	
	\medskip
	
	It turns out that the identifiability of the network, i.e. the ability to recover a module or the whole network from the input-output relation, is a generic notion: Either \emph{almost all} transfer matrices corresponding to a given network structure are identifiable, in which case the structure is called \emph{generically identifiable}, or none of them are. 
	\newb{A number of} works study generic identifiability when all nodes are excited or measured, i.e. when $B$ or $C=I$ \cite{bazanella2017identifiability}, \cite{weerts2018identifiability}. Considering the graph of the network, path-based conditions on the allocation of measurements in the case of full excitation are derived in \cite{hendrickx2018identifiability}. Reformulating these conditions by means of disjoint trees in the graph, \cite{cheng2019allocation} develops a scalable algorithm to allocate excitations/measurements in case of full measurement/excitation.  Rather than verifying the identifiability of a given excitation allocation, \cite{shi2019excitation} addresses the question of where to allocate excitation signals under full measurement.
	
	\medskip
	
	The conditions of \cite{hendrickx2018identifiability} apply for generic identifiability, i.e. identifiability of almost all transfer matrices corresponding to a given network structure.  \cite{van2018topological} extends the path-based conditions under full excitation to determine the identifiability \emph{for all} (nonzero) transfer matrices corresponding to a given structure, and \cite{van2019necessary} provides conditions for the outgoing edges of a node, and the whole network under the same conditions.
	
	\medskip
	
	In all these works, the common assumption \newg{is} that all nodes \newg{are} either excited, or measured. In \cite{bazanella2019network}, this assumption \newg{is} relaxed and generic identifiability with partial excitation and measurement \newg{is} addressed. The authors derive necessary conditions and sufficient conditions for identifying the outgoing and incoming edges of a node, and for particular network topologies such as loops and trees.
	
	\medskip
	
	The question of generic identifiability under partial excitation and measurement for the general case remains unsolved: We wish to find a combinatorial characterization for generic identifiability of an edge or a network, that is expressed purely in terms of graph theoretical properties, akin to what was done in the full excitation case e.g. in \cite{hendrickx2018identifiability}. Such characterization would in particular \newg{pave} the way for optimizing the selection of nodes to be excited and measured, akin to the work in \cite{shi2019excitation} in the full measurement case. 
	\pagebreak
	
	In this paper, we introduce the notion of local identifiability, i.e. we address the question of identifiability only on a neighborhood of $G(z)$. Local identifiability is a necessary condition for identifiability, and we believe its understanding will allow making progress towards understanding identifiability.
	We \newg{show that it is generic,} provide an algebraic necessary and sufficient condition for local identifiability of a network in terms of the rank of a matrix, and a necessary and sufficient condition for local identifiability of an edge depending on the orthogonality of a matrix kernel to \newg{a} certain vector (which is itself equivalent to some rank conditions). These conditions generalize previous results obtained under full excitation or measurement \cite{hendrickx2018identifiability}. 
	We then exploit these conditions to develop an algorithm that determines which transfer functions are locally identifiable, with probability 1. An implementation is available at \cite{matlab}, and includes a graphical representation of the results.
	
	

	
	\medskip
	\noindent\textbf{Assumptions}:
	We consider the problem modeled in \eqref{eq:networkModel}. Consistently with previous works, we assume the invertibility of \newg{the} matrix $(I-G(z))$, \newb{which is equivalent to the network being well-posed}, and that $CT(z)B=C(I-G(z))^{-1}B$ has been identified exactly. We do not suppose having access to any information related to the effect of the noise signals $v_1,v_2$. The additional information that could be gathered from this knowledge in our context is left for further works. 

    \medskip
	
	For the simplicity of exposition, we consider in this paper a single frequency $z$, so that all transfer functions are modeled simply by a complex value, and the matrices $G$ and $T= (I-G)^{-1}$ are complex matrices rather than matrices of transfer functions. Conceptually, our generic results directly extend to the transfer function case: if one can recover a $G_{ij}(z)$ at a given frequency $z$ for almost all $G$ consistent with a network, then one can also recover it at all other frequencies, and hence recover the transfer function. We intend to remove this simplification or to formalize this intuitive argument in a further version of this work. 
	In the remainder of this document, we omit $(z)$ to lighten notations.

	\medskip

	\section{LOCAL IDENTIFIABILITY} \label{sec:local_identif}
	
	\medskip
	
	We \newg{start} by reminding the definition of identifiability, see e.g. {\cite{hendrickx2018identifiability, weerts2018single}}. A network matrix $G$ defines a directed graph where the edge $j\rightarrow i$ is present if $G_{ij} \neq 0$. We say that a matrix $\tilde G$ is consistent with a graph defined by a network matrix $G$ if $\tilde G_{ij}(z)$ is zero when there is no edge $j\rightarrow i$ in the graph defined by $G$.
	
	\medskip
	
		\begin{defn}
			The transfer function \newo{$G_{ij}$} is \emph{identifiable} at $G$ from excitations $\mathcal{B}$ and measurements $\mathcal{C}$ if, for all network matrix $\tilde{G}$ consistent with the graph, there holds
			\begin{align} \label{eq:def1}
				C(I-\tilde{G})^{-1}B = C(I-G)^{-1}B \Rightarrow \tilde{G}_{ij} = G_{ij}.
			\end{align}
			The network matrix is identifiable at $G$ if each transfer function $G_{ij}$ is identifiable at $G$, i.e. if the left-hand side \newg{of {\eqref{eq:def1}}} implies $\tilde G = G$.
			
		\end{defn}
		
		\medskip

	These definitions were made for matrices of transfer functions $G(z)$, but we remind that here we consider a single frequency for simplicity, so $G$ is just a matrix in $\mathbb{C}^{n\times n}$.
	
	\medskip
	
	The combinatorial conditions for identifiability  obtained under full excitation or measurement (i.e. when $B$ or $C = I$) in \cite{gevers2017identifiability}, \newg{{\cite{bazanella2017identifiability}}}, \cite{weerts2018identifiability}, \cite{hendrickx2018identifiability}, \cite{cheng2019allocation} were made possible by a linear algebraic reformulation of the identifiability questions. No equivalent reformulation is available for settings with partial measurements and excitation. Hence we will work on a weaker notion amenable to linear analysis: local identifiability, which corresponds to identifiability provided that $\tilde G$ is sufficiently close to $G$. 
	
	\medskip
	
		\begin{defn}
		The transfer function \newg{$G_{ij}$} is \emph{locally identifiable} \newg{at $G$} from excitations $\mathcal{B}$ and measurements $\mathcal{C}$ if there exists $\epsilon > 0$ such that for any $\tilde{G}$ consistent with the graph satisfying $||\tilde{G}-G||<\epsilon$, there holds
		\begin{align} \label{eq:defLocalIdentifEdge}
			C(I-\tilde{G})^{-1}B= C(I-G)^{-1}B \Rightarrow \tilde{G}_{ij} =G_{ij}.
		\end{align}
		The network matrix is locally identifiable \newb{at $G$} if each transfer function $G_{ij}$ is locally identifiable at $G$, i.e. if the left-hand side of \eqref{eq:defLocalIdentifEdge} implies $\tilde G = G$.
        \end{defn}
	
	\medskip
	
	  Local identifiability is a necessary condition for identifiability, \newb{both for recovering the entire network or just a particular edge}.  It is \emph{a priori} a weaker notion, but we have so far found no example of network that is locally identifiable but not globally identifiable. 
	  Moreover, one can show that if a network is locally identifiable, then it can be recovered up to a discrete ambiguity, i.e. the set of unknown $G$ corresponding to a measured $C(I-G)^{-1}B$ would be discrete. Finally, local identifiability is \newb{relevant} in the many practical situations where the transfer functions are already approximately known.
	  
	  \medskip
	  
	  \noindent\textbf{Genericity:}	Given a graph and sets $\mathcal{B},$ $\mathcal{C}$ of excited and measured nodes, we say that an edge is \emph{generically} (locally) identifiable if it is (locally) identifiable at all $G$ consistent with the graph, except possibly those lying on a \newb{lower-dimensional} set \newo{(i.e. a set of dimension lower than $n$)}. In the remainder of this paper, we will say that a property holds \emph{generically}, or \emph{for almost all} (resp. no) \newb{variables} if it holds for all \newo{(resp. no)} \newb{variables}, except possibly those lying on a lower-dimensional set. 
	  \newg{We will see} that local identifiability is indeed a generic notion: for given network topology and sets of excited and measured nodes, the transfer function corresponding to the edge $(i,j)$ is either identifiable at almost all $G$ (generically identifiable) or identifiable at almost none of them. More details about this will be given in future versions of this work.

	\medskip
	
	\section{IDENTIFIABILITY AND INJECTIVITY} \label{sec:identif_injectiv}
	
	\medskip
	
	The identifiability questions we consider can be reformulated in terms of (local) injectivity. Let $x \in \mathbb{C}^{\abs{E}}$ be a vector compiling the $|E|$ (potentially) nonzero entries of $G$, and $G(x)$ the \newg{network} matrix $G$ corresponding to a vector $x$, so that every matrix consistent with the graph can be written as $G(x)$ for some $x$. 
	We define
	\begin{align} \label{eq:def_f}
			\myfct :\mathbb{C}^{|E|} \setminus \mathcal{D} \rightarrow \mathbb{C}^{|\mathcal{B}|\cdot|\mathcal{C}|}: x \rightarrow \myvec \prt{C(I-G(x))^{-1}B}
	\end{align}
	\new{where the set $\mathcal{D} \triangleq \{ x | \det(I-G(x))= 0 \}$ collects the vectors for which the inverse does not exist.}
	The $|\mathcal{B}|\times|\mathcal{C}|$ image is vectorized into a single $|\mathcal{B}|\cdot|\mathcal{C}|$ vector, which will ease the differential analysis below.
	
	\medskip
	
	Recovering a matrix $G$ from $C(I-G)^{-1}B$ can then be seen as recovering $x$ from $\myfct (x)$, while recovering $G_{ij}$ corresponds to recovering a specific entry $x_e$ of $x$ from $\myfct (x)$ \newg{(where $e$ is the index of $G_{ij}$ in $x$).} 
	\newg{\emph{Local identifiability becomes then local injectivity}}: determining if $\myfct (\tilde x) = \myfct (x)$ implies $\tilde x = x$ for every $\tilde x$ in a neighborhood of $x$. The local identifiability of an edge $e$ becomes a question of local coordinate-injectivity: determining if $\myfct (\tilde x) = \myfct (x)$ implies $\tilde x_e = x_e $ for every $\tilde x$ \newg{in} a neighborhood of $x$. The notion of local coordinate-injectivity is formally introduced in the definition below.
	
	\medskip
	
	\begin{defn} \label{def:e_injectivity}
		The function $f:M \rightarrow N$ is \emph{locally coordinate-injective} for \new{the coordinate} $e$ \newg{at} $x$ if there exists $\epsilon > 0$ such that for all $\tilde{x} \in B(x,\epsilon)$, there holds
		\begin{align} \label{eq:def_e_injectivity}
			f(\tilde{x}) = f(x) \Rightarrow \tilde{x}_e = x_e.
		\end{align}
	\end{defn}
	
	\medskip
	
	
	Intuitively, these local questions could be determined by the gradient: $\myfct$ would be locally injective at $x$ if $\nabla \myfct (x)\cdot \delta \neq 0$ for any nonzero $\delta \in \mathbb{C}^{\abs{E}}$, i.e. if $\nabla \myfct (x)$ has a trivial kernel. Observe this requires the number $|\mathcal{B}|\cdot|\mathcal{C}|$ of available input-output relations to be at least as large as the number $\abs{E}$ of unknowns. Similarly, $\myfct$ would be coordinate-injective at $x$ for edge $e$ if $\nabla \myfct (x)\cdot \delta \neq 0$ for any $\delta \in \mathbb{C}^{\abs{E}}$ with $\delta_e \neq 0$, i.e. if the kernel of $\nabla \myfct (x)$ is orthogonal to $\evec_e$, \newb{the standard basis vector filled with zeros except 1 at the $e$-th entry.}
	
	\medskip
	
	However, the issue is more complex than it first seems. For example, the real function $g(x) = x^3$ is injective, but $g'(0) = 0$ admits a non-trivial kernel, in contradiction with our previous statement. The problem is that the derivative of $g$ vanishes at $0$, while being nonzero everywhere else. Adding an assumption on the conservation of information contained in the derivative (i.e. the rank of its gradient) resolves that issue.
	The next lemma, \newg{proved in the Appendix,} shows that the intuitive idea behind the gradient is in fact valid when $\nabla f$ has constant rank. We will then see in Proposition \ref{prop:kerdiff_perp_ee} how this rather constraining assumption can be removed provided one is only interested in local injectivity \emph{at almost all $x$}.
    
    \medskip
    
	\begin{lemma} \label{lemma:kerdiff_perp_ee}
		Suppose $f: M \rightarrow N$ is $C^\infty$ where \new{$M$ and $N$ are smooth manifolds, $M$ or $N$ has finite dimension}, and $\nabla f(x)$ has the same rank for every $x \in M$. 
		Then $f$ is locally coordinate-injective for $e$ \newg{at all $x \in M$} if and only if
		for all $x \in M$,
		\begin{align} \label{eq:kerdiff_perp_ee}
			\ker \nabla f(x) \perp \evec_e,
		\end{align}
		\newb{where $\evec_e$ denotes the standard basis vector filled with zeros except 1 at the $e$-th entry.}
	\end{lemma}
	
	\medskip
	
	The next proposition, \newg{proved in the Appendix,} characterizes \newg{local} injectivity almost everywhere for analytic functions.  Analytic functions are functions which can be locally approximated everywhere by their Taylor expansion, see e.g. \new{{\cite{hille2005analytic}}} for a precise definition. They include in particular all compositions of products, sums, divisions, etc.
	
	\medskip

		\begin{prop} \label{prop:kerdiff_perp_ee}
		Let $f: M \rightarrow N$ be an analytic function where \new{$M$ and $N$ are smooth manifolds, $M$ or $N$ has finite dimension, and $M$ is open}. \newb{Then exactly one of the two following holds:}
		\begin{enumerate}[(i)]
		    \item \newb{$\ker \nabla f(x) \perp \evec_e$ for almost \emph{all} $x$ and $f$ is locally coordinate-injective for $e$} \newg{at almost \emph{all} $x$};
		    \item \newb{$\ker \nabla f(x) \perp \evec_e$ for almost \emph{no} $x$ and $f$ is locally coordinate-injective for $e$ at almost \emph{no} $x$.}
		\end{enumerate}

    \end{prop}
    
    \medskip
	
	\section{NECESSARY AND SUFFICIENT CONDITIONS FOR LOCAL IDENTIFIABILITY} \label{sec:conditions}

    \medskip
    
    \newg{First we define the $\abs{\mathcal{B}}\abs{\mathcal{C}}\times |E|$ matrix $K$ which, as will be proved in Lemma {\ref{lemma:differential_f}}, is the gradient of $\myfct$.}
        \begin{align} \label{eq:def_K}
			K(x) \triangleq \prt{B^T T^T(x) \otimes CT(x)} I_G,
		\end{align}
		where $T(x) \triangleq (I-G(x))^{-1}$, $\otimes$ denotes the Kronecker product and the matrix $I_G \in \{0,1\}^{n^2 \times |E|}$ selects only the columns of the preceding \new{$\abs{\mathcal{B}}\abs{\mathcal{C}}$}$\times n^2$ matrix corresponding to actual edges. It is defined by
		\begin{align*}
			I_G \triangleq \left[ \myvec(G(\evec_1)) \ \myvec(G(\evec_2)) \, \cdots \, \myvec(G(\evec_{|E|})) \right].
		\end{align*}

    \medskip

    We are now ready to \newg{show that local identifiability is generic, and}
    derive our necessary and sufficient condition. 

    \bigskip
    
	\begin{thm}\label{thm:main_alg_condition}
	    \newb{Let $K$ be as defined in {\eqref{eq:def_K}} and $e$ be the index of $G_{ij}$ in $x$. Exactly one of the two following holds:}
	    \begin{enumerate}[(i)]
		    \item \newb{$\ker K(x) \perp \evec_e$ for almost \emph{all} $x$ and $G_{ij}$ is generically locally identifiable;}
		    \item \newb{$\ker K(x) \perp \evec_e$ for almost \emph{no} $x$ and $G_{ij}$ is generically locally non-identifiable,} \newg{hence generically non-identifiable.}
		\end{enumerate}
		\newo{Moreover, $\ker K(x) \perp \evec_e$ is equivalent to the following implication holding for all $\Delta$ consistent with the graph:}
		\begin{align} \label{eq:CTDeltaTB_Delta}
			C (I-G)^{-1} \Delta (I-G)^{-1} B = 0 \Rightarrow \Delta_{ij} = 0.
		\end{align}
		
	\end{thm}
	
	\bigskip


    \newg{The proof of Theorem {\ref{thm:main_alg_condition}} relies on applying Proposition {\ref{prop:kerdiff_perp_ee}} on the gradient of $\myfct$, characterized by Lemma {\ref{lemma:differential_f}}.}

    \bigskip

	\begin{lemma} \label{lemma:differential_f}
		Let $\myfct$ be the function defined in \eqref{eq:def_f}. Its gradient is given by
		\begin{align} \label{eq:grad_f}
			\nabla \myfct (x) = \prt{B^T T^T(x) \otimes CT(x)} I_G = K(x).
		\end{align}
		Besides, for all $\delta \in \mathbb{C}^{|E|}$, there holds
		\begin{align} \label{eq:differential_f}
			\mymat \prt{\nabla \myfct (x) \, \delta}
			= C T(x) \, G(\delta) \, T(x) B,
		\end{align}
		where $\mymat$ reorganizes its $|\mathcal{B}||\mathcal{C}|$ vector in a $|\mathcal{B}|\times|\mathcal{C}|$ matrix.
		\end{lemma}
	
	\pagebreak
	
	\begin{proof}
	First note that $\nabla \myfct (x) \, \delta$ is the differential of $\myfct$ at $x$ evaluated on $\delta$, which we denote $\diff_x \myfct (\delta)$.
	    We start by deriving the expression of the differential \eqref{eq:differential_f}, from which we compute the gradient \eqref{eq:grad_f}.
	    \newg{The differential of $\myfct$ follows from the differential of $T$ by linearity:} $\mymat (\diff \myfct) = C \diff T B,$
	    \newo{which relies on its partial derivatives, developed as {\cite{petersen2012matrix}}:} 
		\begin{align*}
			\fpart{T}{x_e}(x)
			&= - \prt{I-G(x)}^{-1} \fpart{(I-G)}{x_e}(x) \prt{I-G(x)}^{-1}\\
			&= - T(x) \prt{-\fpart{G}{x_e}(x)} T(x) = T(x) \, G(\evec_e) \, T(x).
		\end{align*}
		\newg{The partial derivatives then yield the differential:}
		\begin{align*}
			\diff_x T(\delta)
			= \sum_{e=1}^{|E|} T(x) \, G(\evec_e) \, T(x) \, \delta_e
			= T(x) \, G(\delta) \, T(x).
		\end{align*}
		For the differential of $\myfct$, we obtain
		\begin{align} \label{eq:differential_f_2}
			\mymat (\diff_x \myfct (\delta))
			= C \diff_x T(\delta) B
			= C T(x) \, G(\delta) \, T(x) B.
		\end{align}
		In order to derive the gradient of $\myfct$, we vectorize \eqref{eq:differential_f_2}:
		\begin{align*}
			\nabla \myfct (x) \, \delta
			&= \myvec \prt{ C T(x) \, G(\delta) \, T(x) B} \\
			&= \prt{B^T T^T(x) \otimes CT(x)} \myvec \prt{G(\delta)},
		\end{align*}
		where $\myvec$ stacks the columns of its $|\mathcal{B}|\times|\mathcal{C}|$ matrix argument into a $|\mathcal{B}|\cdot|\mathcal{C}|$ vector.
		Then, observing that $\myvec \prt{G(\delta)} = I_G \, \delta$ gives
		\newo{$\nabla \myfct (x) \, \delta = K(x) \, \delta \quad \forall \, \delta \in \mathbb{C}^{|E|}$, which yields {\eqref{eq:grad_f}}.}
    \end{proof}

	\medskip

	\textbf{Proof of Theorem \ref{thm:main_alg_condition}:} 
		\newg{As explained in Section {\ref{sec:identif_injectiv}}, the generic local identifiability of transfer function $G_{ij}$ is equivalent to function $\myfct$ defined in {\eqref{eq:def_f}} being locally coordinate-injective for $e$ at almost all $x$ ($e$ is the index of $G_{ij}$ in $x$).}
		
		\medskip
		
		Proposition \ref{prop:kerdiff_perp_ee} gives conditions on local coordinate-injectivity, let us verify its assumptions \newg{for $\myfct$}.
		\newo{First $\myfct$ is analytic, and its domain and image sets have finite dimension.}
		\newb{The openness of the domain $\mathbb{C^{|E|}}\setminus \mathcal{D}$ follows from the continuity of $\det (I-G)$. Since the inverse image of a closed set by a continuous function is a closed set {\cite{hausdorff1991set}}, the zero set of $\det (I-G)$ is closed, and its complement is open. Its assumptions being fulfilled, Proposition {\ref{prop:kerdiff_perp_ee}} applies. The first part of the theorem then follows directly from {\eqref{eq:grad_f}}.}
		
		\medskip
		
		Besides, $\ker \nabla \myfct (x) \perp \evec_e$ can be rewritten equivalently as
		\begin{align*}
		    \nabla \myfct (x) \, \delta = 0 \Rightarrow \delta_e = 0.
		\end{align*}
		Reorganizing this $|\mathcal{B}|\cdot|\mathcal{C}|$ vector into a $|\mathcal{B}|\times|\mathcal{C}|$ matrix \newg{gives}
		\begin{align*}
		    \mymat \prt{ \nabla \myfct (x) \, \delta } = 0 \Rightarrow \delta_e = 0,
		\end{align*}
		and \newo{combining} it with \eqref{eq:differential_f} then yields \eqref{eq:CTDeltaTB_Delta}.
	\cqfd

	Theorem \ref{thm:main_alg_condition} provides a necessary and sufficient condition for local identifiability of a transfer function. Therefore, combining this condition for all edges immediately leads to a characterization of local identifiability of the whole network.
	
	\pagebreak
	
	\begin{corol} \label{corol:network_identif}
	    \newb{Exactly one of the two following holds:}
	    \begin{enumerate}[(i)]
		    \item \newb{$\rank K(x) = |E|$ for almost all $x$ and $G$ is generically locally identifiable;}
		    \item \newb{$\rank K(x) = |E|$ for almost \emph{no} $x$ and $G$ is generically locally non-identifiable.}
		\end{enumerate}
		\newb{Moreover, $\rank K(x) = |E|$ is equivalent to the following implication holding for all $\Delta$ consistent with the graph:}
		\begin{align} \label{eq:CTDeltaTB_Delta_net}
			C (I-G)^{-1} \Delta (I-G)^{-1} B = 0 \Rightarrow \Delta = 0.
		\end{align}
	\end{corol}
	
	\medskip
		
	\begin{proof}
		\newg{First,} $\ker K(x) \perp \evec_e$ for each $e$ \newg{is equivalent to} $\ker K(x) = \{ 0 \}$, and therefore $\rank K(x) = |E|$ by the rank-nullity theorem.
		\newg{Then, {\eqref{eq:CTDeltaTB_Delta}} for each $(i,j)$ yields {\eqref{eq:CTDeltaTB_Delta_net}}}.
	\end{proof}

	
	We observe that condition \eqref{eq:CTDeltaTB_Delta_net} reduces to the necessary and sufficient condition for (global) identifiability in the full excitation (resp. measurement) case. Indeed, since $T = (I-G)^{-1}$ is by construction invertible, $CT\Delta TB = 0$ is equivalent to $CT\Delta=0$ when $B$ is the identity matrix (resp. to $\Delta T B$ when $C$ is the identity matrix), as in \cite{hendrickx2018identifiability}.

	\medskip
	
	\section{ALGORITHM} \label{sec:algorithm}
	
	\medskip
	
	The necessary and sufficient conditions of Theorem \ref{thm:main_alg_condition} can be tested exactly by symbolic computation, but it becomes rapidly computationally intractable. Hence, we have designed a probability-1 algorithm that exploits the generic character \newb{of local identifiability, as proved in Theorem {\ref{thm:main_alg_condition}}.}
	Testing these conditions with a randomly selected $x\in \mathbb{C}^{\abs{E}}$ will therefore provide the correct result with a probability 1, and can be done very efficiently. \newg{This extends what E. J. Davison does in the controllability context {\cite{davison1977connectability}}.}
	
	\medskip
	
	Although this idea allows to design an algorithm working with probability 1, we cannot exclude the possibility that an actual implementation could suffer from numerical issues for large networks, with e.g. the numerically computed rank being very close to zero if the randomly selected values are close to the problematic lower-dimensional set. Hence we repeat the procedure several times to increase the reliability of our results, as described in Algorithm \ref{algo:identif_test}. However, we did not observe any such numerical problem in any of our tests. Our implementation is available at \cite{matlab}, and includes a graphical representation of the results.
	
	\begin{algorithm}
    \caption{Identifiability test}\label{algo:identif_test}
      \begin{algorithmic}[1]
        \Require Graph topology, matrices $B$ and $C$, \mytexttt{nsamples}
        \Ensure The identifiability of each transfer function
        \State Initialize \mytexttt{network} to \mytexttt{false} and \mytexttt{edges} to a \mytexttt{false} vector of length $|E|$ (\mytexttt{true} means identifiable)
        \For{$i \leftarrow 1$ \textbf{to} \mytexttt{nsamples} \textbf{by} 1}
            \State Randomly generate a complex network matrix $G$
            \State Construct $K= \prt{B^TT^T \otimes CT} I_G$
            \If{$\rank K = |E|$}
                \State \mytexttt{network} $\leftarrow$ \mytexttt{true}
            \Else
                \ Compute $V$, basis of $\ker K$
                \State Calculate $v$, the binary vector of length $|E|$, with $v_e =$ \mytexttt{true}  if the $e$-th entry of each vector of $V$ is 0, $v_e =$ \mytexttt{false} otherwise
                \State \mytexttt{edges} $\leftarrow$ \mytexttt{edges} \textbf{or} $v$ \textbf{entry-wise}
            \EndIf
        \EndFor
        
        \If{\mytexttt{network}} \textbf{return} Network identifiable
        \Else \textbf{ return} \mytexttt{edges}
        \EndIf
      \end{algorithmic}
    \end{algorithm}
    
    \medskip
    
\begin{figure}
    \centering
    \subfloat[]{\label{subfig1}\includegraphics[width=\linewidth]{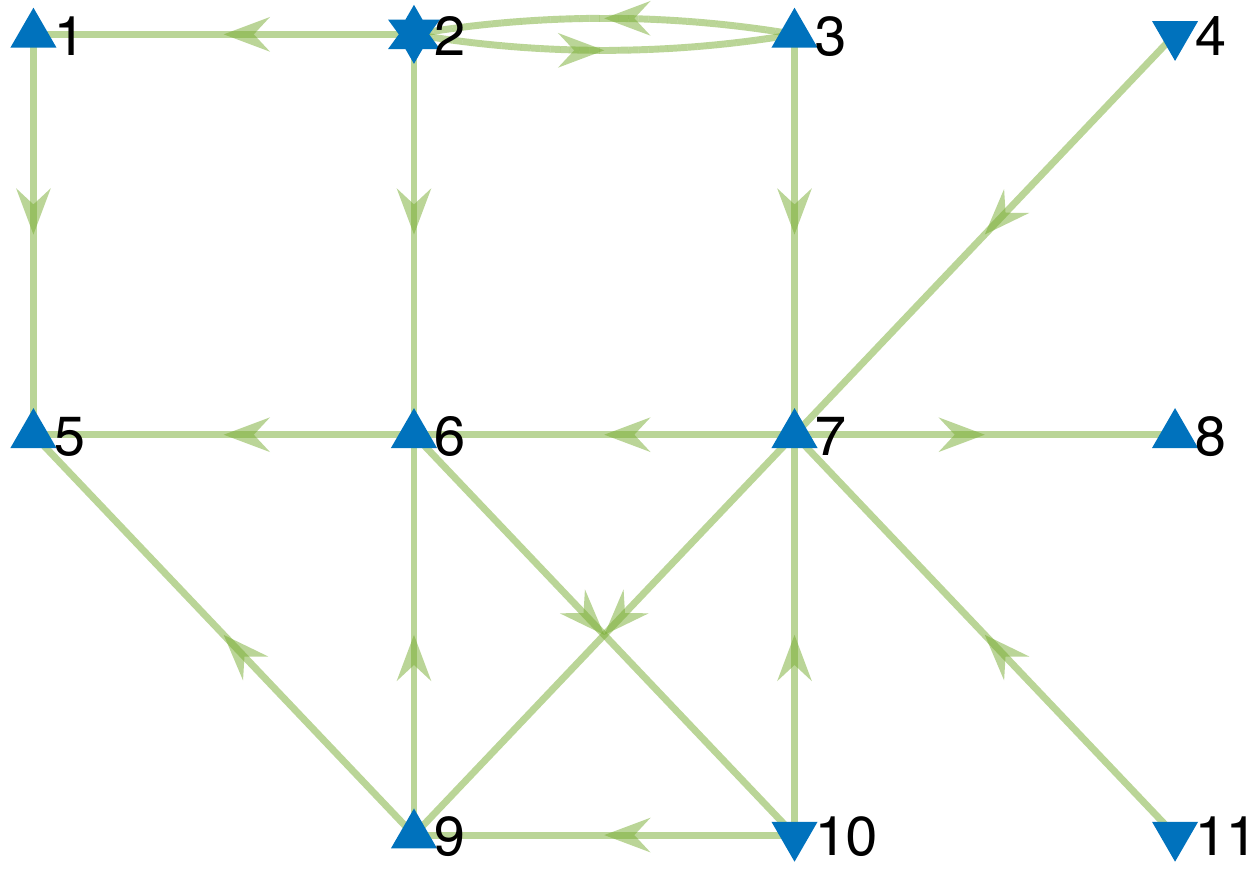}} \\
    \subfloat[]{\label{subfig2}\includegraphics[width=\linewidth]{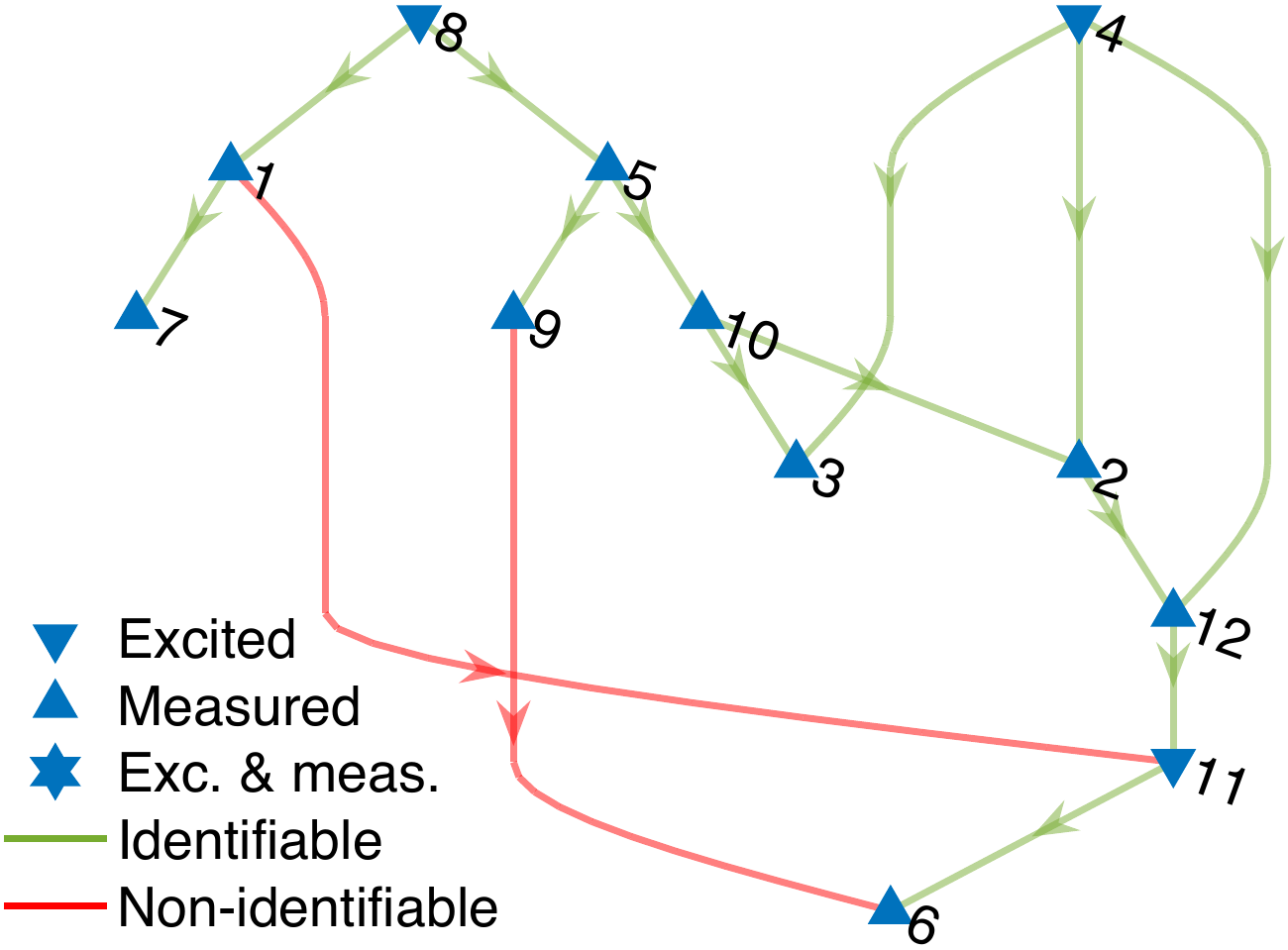}} 
    \caption{Example of outputs from Algorithm \ref{algo:identif_test}: (a) is an improvement on Example \newg{2} of \cite{cheng2019allocation} showing not all nodes need to be measured.  In (b), two edges cannot be identified, but exciting 1 or 9 restores full local identifiability.} 
    \label{fig:tests_matlab}
    \end{figure}

    Figure \ref{fig:tests_matlab} shows examples of outputs of Algorithm \ref{algo:identif_test}.
	Figure \ref{fig:tests_matlab}(a) is the network of Example \newg{2} in \cite{cheng2019allocation}, for which the authors \newg{provide} an optimal excitation set assuming all nodes \newg{are} measured. Our algorithm shows that all edges can still be locally identified if some of the nodes (i.e. 4, 10, 11) are not actually measured, \newg{and node 9 must not be excited.} Note this does not contradict the optimality of the excited set found in \cite{cheng2019allocation}, since this set was found under the assumption of full measurement, \newg{and node 9 was affected by external noise.}
	
	\medskip
	
	Figure \ref{fig:tests_matlab}(b) shows a nontrivial example of network where all edges but two are generically locally identifiable. Interestingly, exciting 1 or 9 allows recovering both (1, 11) and (9, 6), even though these edges do not appear directly related. 
	
	
	\section{CONCLUSION}
	
	This work was motivated by one main open question: determining graph-theoretical or combinatorial conditions for generic network and edge identifiability in networked systems. 
	
	The local notions we have introduced allowed making progress in this direction. In particular, the necessary and sufficient conditions for generic local identifiability in terms of generic ranks allowed us to design an algorithm determining, with probability 1, the generic local identifiability of each edge. 
	
	Providing a graph-theoretical characterization for local identifiability remains an open question, but we believe that our generic rank-based characterization could pave the way to a solution. Moreover, our algorithm allows rapidly testing ideas and conjectures, hence facilitates further research.  

	In principle, generic local identifiability is a weaker notion than generic identifiability. However, we were so far unable to find examples of networks that are locally, but not globally, identifiable. The possible equivalence of the two notions would greatly simplify the study of identifiability, and remains an open question. 

	\section{ACKNOWLEDGMENTS}

	The authors gratefully acknowledge P.-A. Absil for his \newg{help on manifold theory}, and M. Gevers and A. S. Bazanella for the interesting and insightful discussions.

	\bibliographystyle{ieeetr}
	\bibliography{cdc20}
	
	\appendix
	\subsection{Proof of Lemma \ref{lemma:kerdiff_perp_ee}}
	
		We must prove that \eqref{eq:def_e_injectivity} \newb{ for all $x \in M$} $\Leftrightarrow$ \eqref{eq:kerdiff_perp_ee} for all $x \in M$. \new{The smoothness of $f$ and constant rank of $\nabla f$ allow the application of the subimmersion theorem {\cite{abraham2012manifolds}}. In short, this theorem states that the level set of an image $n_0$ by such function is a submanifold of $M$, and its tangent space coincides with the kernel of the gradient. Here, we take $n_0 = f(x)$ to make the connection with injectivity,} and the theorem yields that for all $x$, $f^{-1}(f(x))$ is a submanifold of $M$ with
		\begin{align} \label{eq:subimmersion_general_result}
			T_{\tilde{x}} f^{-1}(f(x)) = \ker \nabla f(\tilde{x}) \quad \forall \ \tilde{x} \in f^{-1}(f(x)),
		\end{align}
		where $T_{\tilde{x}} A$ denotes the tangent space to manifold $A$ at $\tilde{x}$.
		
		\medskip
		
		We start with necessity: \newb{assume {\eqref{eq:def_e_injectivity}} for all $x$, and consider a specific $x$}. Each $\tilde{x}$ belonging to a neighborhood of $x$ in  $f^{-1}(f(x))$ has the same $e$-coordinate as $x$. Therefore, any tangent vector to $f^{-1}(f(x))$ at $x$ has a zero $e$-component, and the tangent space to this manifold at $x$ is orthogonal to direction $e$, i.e. $T_x f^{-1}(f(x)) \perp \evec_e$. Equation \eqref{eq:subimmersion_general_result} evaluated at $\tilde{x} = x$ then gives \eqref{eq:kerdiff_perp_ee} and proves necessity.
		
		\medskip

		We now prove sufficiency: \newb{assume {\eqref{eq:kerdiff_perp_ee}} for all $x$, and consider a specific $x$.} Combining \eqref{eq:kerdiff_perp_ee} with \eqref{eq:subimmersion_general_result} yields $T_{\tilde{x}} f^{-1}(f(x)) \perp \evec_e$ for all $\tilde{x} \in f^{-1}(f(x))$.
		\newg{Take a specific $\tilde{x} \in f^{-1}(f(x))$, and assume that there is a path lying in $f^{-1}(f(x))$ between $x$ and $\tilde x$.}
		At every point of this path, the tangent space is orthogonal to direction $e$, so the derivative along the path is zero in the $e$-coordinate. Integrating along the path gives a zero increase in \newo{$e$-coordinate, hence $\tilde{x}_e = x_e$.}
		
		\medskip
		
		\newg{To complete the proof, we need to show that there exists $\epsilon > 0$ such that for all $\tilde x \in B(x,\epsilon) \cap \ f^{-1}(f(x))$, there is a path lying in $f^{-1}(f(x))$ between $x$ and $\tilde x$.} \newo{Such path can be found with the rank theorem {\cite{abraham2012manifolds}}}.
		\newb{Roughly speaking, this theorem states that if $\nabla f$ has constant rank $k$ in a manifold of dimension $m$, then $m-k$ variables are redundant and can be eliminated.}
		\neww{More precisely, for all $p$ in the manifold, there exist local coordinates centered at $p$ and $f(p)$ in which $f(q)=(q_1, \dots, q_k, 0, \dots, 0)$ for any $q$ in the manifold.
		
		\medskip
		
		For the rest of this proof, we work in the local coordinates provided by this theorem, centered at $p=x$. Since these are centered at $x$ and $f(x)$, $x_1=\dots=x_k=0.$}
		\newg{Let us now take $\epsilon > 0$, and $\tilde x \in B(x,\epsilon) \cap \ f^{-1}(f(x))$.}
		\neww{The rank theorem yields that the last entries of $f(\tilde x)$ are zero in those local coordinates, i.e. $f(\tilde x) = (\tilde x_1, \dots, \tilde x_k, 0, \dots, 0)$. Since $\tilde x$ lies in the manifold $f^{-1}(f(x))$, it must have the same image as $x$, hence $f(\tilde x) = f(x) = 0$, and $\tilde x_1 = \dots = \tilde x_k = 0$.
		An example of path between $x$ and $\tilde x$ is given by $(0, \dots, 0, \lambda \tilde x_{k+1}, \dots, \lambda \tilde x_{|E|}$), with $\lambda \in [0, 1]$. Applying the function $f$ on any point of this path gives $0$ in the local coordinates by the rank theorem, showing that the path is indeed included in the manifold $f^{-1}(f(x))$.}
		 \cqfd
		 

    \subsection{Proof of Proposition \ref{prop:kerdiff_perp_ee}}
    
    \newg{The proof of Proposition {\ref{prop:kerdiff_perp_ee}} relies on Lemma {\ref{lemma:cstt_rank_ker_grad}}.}
    
    \begin{lemmaappendix} \label{lemma:cstt_rank_ker_grad}
        \newg{The gradient of an analytic function $f$ reaches its maximal rank everywhere except on a closed lower-dimensional set.}
	    \newb{Moreover, the orthogonality relation}
	    \begin{align} \label{eq:kerdiff_perp_ee_2}
			\ker \nabla f(x) \perp \evec_e
		\end{align}
		\newb{either holds:}
		\begin{enumerate}[(i)]
		    \item \newb{for all $x$ except those lying on a closed lower-d. set;}
		    \item \newb{only for $x$ lying on a closed lower-dimensional set.}
		\end{enumerate}
    \end{lemmaappendix}
    
    \medskip
    
    \textbf{Sketch of proof:}
    \newg{First, the set on which the rank of a matrix drops can be expressed as the intersection of zero sets of determinants of submatrices {\cite{van1991graph}}. The determinant is analytic in its entries, which are analytic in $x$ since $f$ is assumed analytic (hence so is $\nabla f$). Since nonconstant analytic functions vanish only on a lower-dimensional set {\cite{d2010introduction}}, the set on which $\nabla f(x)$ does not reach its maximal rank has lower dimension. Besides, this set is also closed: the inverse image of a closed set by a continuous function is a closed set {\cite{hausdorff1991set}}, hence the zero sets of the determinants are closed, and so is their intersection. This extends the notion of generic rank introduced in {\cite{davison1977connectability}}.}

    \newg{Then, {\eqref{eq:kerdiff_perp_ee_2}} can be equivalently rewritten as $\nabla f(x) \delta = 0 \Rightarrow \delta_e = 0 \ \forall \ \delta$, which is equivalent to the linear independence of the $e$-th column of $\nabla f(x)$ with the other columns. It can be formulated by means of ranks, and a reasoning similar to above completes the proof.}	\cqfd
    

	
	\textbf{Proof of Proposition \ref{prop:kerdiff_perp_ee}:}
	\newg{We denote $\tilde M$ the intersection of the set on which $\nabla f$ reaches its maximal rank (hence constant), and the set on which {\eqref{eq:kerdiff_perp_ee_2}} holds (case (i)) or does not hold (case (ii)). Lemma {\ref{lemma:cstt_rank_ker_grad}} ensures that those two sets cover the whole domain $M$ except a closed lower-dimensional set. Hence, so does $\tilde M$, and it is open since $M$ is assumed open. Thus, there is an $\epsilon$ for every $x \in \tilde M$ such that the rank of $\nabla f$ is constant over $B(x,\epsilon)$, and {\eqref{eq:kerdiff_perp_ee_2}} holds (case (i)) or does not hold (case (ii)) on the ball. Denote $f_{x,\epsilon}$ the restriction of $f$ on $B(x,\epsilon)$.
	The result follows by applying Lemma {\ref{lemma:kerdiff_perp_ee}} on $f_{x,\epsilon}$ for all $x\in \tilde M$.} \hfill \rule{2mm}{2mm}

\end{document}